\documentstyle[amstex,epsf]{article}

\begin{document}

\baselineskip = 15 pt
\parskip 4pt

\newtheorem{prop}{Proposition}
\newtheorem{cor}{Corollary}
\newtheorem{defin}{Definition}
\newtheorem{rem}{Remark}
\newcommand{\Al}{${\cal A}_{\hbar,\eta}(\widehat{sl}_2)~$}
\newcommand{\Alo}{${\cal A}_{\hbar,0}(\widehat{g})~$}
\newcommand{\Apq}{${\cal A}_{q,p}(\widehat{sl}_2)~$}
\newcommand{\Apqg}{${\cal A}_{q,p}(\widehat{g})~$}
\newcommand{\Alg}{${\cal A}_{\hbar,\eta}(\widehat{g})~$}
\newcommand{\Algi}[2]{{\cal A}_{\hbar,#1}(\widehat{g})_{#2}}
\newcommand{\AlN}{${\cal A}_{\hbar,\eta}(\widehat{sl}_N)~$}

\rightline{IMP-NWU/970325}
\rightline{q-alg/9703046}
\rightline{Revised June 6, 1997}

\centerline{\LARGE \bf The algebra ${\cal A}_{\hbar,\eta}(\widehat{g})$ }
\centerline{\LARGE \bf
and Infinite Hopf family of algebras}

\vspace{1cm}
\centerline{Bo-Yu Hou \hspace{1cm}  Liu Zhao}
\centerline{Institute of Modern Physics, Northwest University, Xian 710069,
China}
\vspace{0.3cm}
\centerline{Xiang-Mao Ding}
\centerline{Institute of Theoretical Physics,
Academy of China, Beijing 100080, China}

\date{}

\vspace{3cm}

\begin{abstract}
New deformed affine algebras \Alg are defined
for any simply-laced classical Lie algebra $g$, which are
generalizations of the algebra \Al recently proposed by
Khoroshkin, Lebedev and Pakuliak (KLP). Unlike the work of
KLP, we associate to the new algebras the structure of an infinite
Hopf family of algebras in contrast to the one containing
only finite number of algebras introduced by KLP.
Bosonic representation for \Alg at level 1 is obtained, and
it is shown that, by repeated application of Drinfeld-like
comultiplications, a realization of \Alg at any positive
integer level can be obtained. For the special case of $g=sl_{r+1}$,
$(r+1)$-dimensional evaluation representation is given. The
corresponding intertwining operators are defined and the
intertwining relations are also derived explicitly.
\end{abstract}

\newpage

\section{Introduction}
Since Drinfeld \cite{D1,D2,D3} proposed the quantum groups and Yangian
algebras as deformations of the universal envoloping
algebras of the classical Lie algebras, Hopf algebras with
nontrivial coalgebra structure, especially $q$-affine algebras \cite{RS}
and Yangian doubles \cite{KT,K}, have become one of the
major subjects of pure and applied algebra studies.
Recent progress in the study of Hopf algebras and applications
include the free boson representations of $q$-affine algebras
and Yangian doubles at higher level \cite{slN,HZD,KNO} and the
possibility of describing the dynamical symmetries and
solving the correlation functions
of certain solvable lattice statistic models and integrable
quantum field theories within a pure algebraic framework
\cite{AJMP,BL,Jimbo1,KNO2,Lk1,LS,LuPu,S2}.
The latter problem is, if not the sole force, among the driving
forces which lead to the studies of deformed algebras beyond Hopf
algebras. Examples of such deformed algebras are $q$-
\cite{AKOS,AKOS2,FF,FR,LUKU,LuPu,eg,SKAO,SKMAO} and
$\hbar$- \cite{DHZ,YWL} deformed Virasoro and $W$ algebras,
the elliptic algebra
\Apq \cite{Foda1,Foda2} and its scaling limit \Al \cite{KLP}, and the
algebra of screening
operators of the $q$-deformed $W$-algebras \cite{FJM}
and so on.

In this paper we extend the recent work of
S.Khoroshkin, D.Lebedev and S.Pakuliak \cite{KLP} on
the scaling algebra \Al of the elliptic algebra
\Apq to the general case, \Alg, where $g$ can be any
classical simply-laced Lie algebra of any admissible rank. The
algebra \Al introduced in \cite{KLP} is a formal algebra
with generators carrying {\it continuous} indices.
One of the principal motivation of \cite{KLP}
was to establish a better understanding of the algebra
\Apq from the representation theoretic point of view
because the representation theory of \Apq has been rather
obscure since its birth \cite{Foda1,Foda2}. For this the authors
of \cite{KLP} considered the scaling limit, \Al, instead of \Apq itself,
with generating functions being analytic along some
strip--which plays the role of fundamental parallelogram
for the elliptic algebra \Apq--in the complex plane.
The algebra \Al turns out to be {\it not}
a Hopf algebra but belongs to a Hopf family of algebras in which
the comultiplication can be made associative but with the sacrifice
of changing the periods of structure functions for different iterations
of the comultiplication.
Moreover the twisted intertwining operators appeared in the
representation theory of the algebra \Al satisfy a familar
set of commutation relations which were used in the calculation
of correlation functions for Sine-Gordon model.

We shall show that the algebra \Al actually belongs to (and
constitutes the simplest example of) a new type
of deformed affine algebras, \Alg. Just like their simplest
representative, \Al, these new deformed affine
algebras are not Hopf algebras, because the second deformation
parameter $\eta$ spoils the usual Hopf algebra structure.
However, for two reasons we regard them as {\em deformations} of the
usual Hopf algebras. First, as the second deformation parameter
$\eta$ approaches zero, the currents for the algebra
\Alg obey the same commutation relations as that of the
Yangian double with center $DY_\hbar(g)_c$; Second,
if we consider the level zero representation of \Alg, the
algebraic relations become Hopf algebra again.

Due to the complication for the case of general $g$, we pertain ourselves
only in the current realization. In this form, it is not easy
to write down the analog of comultiplication used by KLP \cite{KLP}.
We therefore use an analog of the well-known Drinfeld comultiplication
to study some aspects in the structure of our algebra. It is, however,
not known whether the finite Hopf family structure of KLP can be
realizad using this form of comultiplication. To cercomevent this drawback,
we introduce an alternative notion, which we call the
infinite Hopf family of algebras, to write down the interations of
comultiplications in a convenient form. It turns out that this
new notion leads to an astonishing simple realization of our algebra
at any positive integer level.

Besides the
pure algebraic elegance, our algebras are also
expected to have relevant applications in such fields as
to develop an algebraic formulation of quantum symmetry
and the calculation of correlation functions for
affine Toda field theories.

\section{The algebras \Alg and infinite Hopf families}

\subsection{The algebra \Alg}

We begin our study by introducing the formal current algebra
(denoted also by the symbol \Alg) associated with \Alg. The
special case of $g=sl_2$ can be inferred from \cite{KLP}.
For other simply-laced classical Lie algebras $g$, the following
definition is, to our knowledge, not introduced anywhere else.

\begin{defin}
The current algebra \Alg associated with the classical simply-laced
Lie algebra $g$ of rank $r$ (as an associative algebra with unit
over the field $\Bbb{C}$) is generated by the $3r$ currents
$\{H^{\pm}_i(u),~E_i(u),~F_i(u)| i=1,~...,~r \}$, the
central element $c$ and 1 with the following generating relations
\footnote{Throughout this paper, the suffices $i$ of the current
operators take integer values, which
indicate different root directions of the underlying Lie
algebra $g$, whilst the symbol
$i$ preceding the $\hbar$'s in the structure functions are
square root of $-1$.}:
\end{defin}

\begin{eqnarray}
H^+_i(u)H^-_j(v)
\!\!\!\!&=&\!\!\!\! \frac{\mbox{sh}\pi\eta(u-v-i\hbar(B_{ij}-c/2))}
{\mbox{sh}\pi\eta(u-v+i\hbar(B_{ij}+c/2))}
\frac{\mbox{sh}\pi\eta'(u-v+i\hbar(B_{ij}-c/2))}
{\mbox{sh}\pi\eta'(u-v-i\hbar(B_{ij}+c/2))}
H^-_j(v)H^+_i(u), \label{CurrB} \\
H^\pm_i(u)H^\pm_j(v)
\!\!\!\!&=&\!\!\!\! \frac{\mbox{sh}\pi\eta(u-v-i\hbar B_{ij})}
{\mbox{sh}\pi\eta(u-v+i\hbar B_{ij})}
\frac{\mbox{sh}\pi\eta'(u-v+i\hbar B_{ij})}
{\mbox{sh}\pi\eta'(u-v-i\hbar B_{ij})}
H^\pm_j(v)H^\pm_i(u),\\
H^\pm_i(u) E_j(v)
\!\!\!\!&=&\!\!\!\! \frac{\mbox{sh}\pi\eta(u-v-i\hbar(B_{ij} \mp c/4))}
{\mbox{sh}\pi\eta(u-v+i\hbar(B_{ij} \pm c/4))}
E_j(v)H^\pm_i(u),\\
H^\pm_i(u) F_j(v)
\!\!\!\!&=&\!\!\!\! \frac{\mbox{sh}\pi\eta'(u-v-i\hbar(B_{ij} \mp c/4))}
{\mbox{sh}\pi\eta'(u-v+i\hbar(B_{ij} \pm c/4))}
F_j(v)H^\pm_i(u),
\end{eqnarray}
\begin{eqnarray}
E_i(u) E_j(v)
\!\!\!\!&=&\!\!\!\! \frac{\mbox{sh}\pi\eta(u-v-i\hbar B_{ij})}
{\mbox{sh}\pi\eta(u-v+i\hbar B_{ij})}
E_j(v)E_i(u),\\
F_i(u) F_j(v)
\!\!\!\!&=&\!\!\!\! \frac{\mbox{sh}\pi\eta'(u-v + i\hbar B_{ij})}
{\mbox{sh}\pi\eta'(u-v - i\hbar B_{ij})}
F_j(v)F_i(u),\\
{}[ E_i(u),F_j(v) ]
\!\!\!\!&=&\!\!\!\! \frac{2\pi}{\hbar} \delta_{ij} \left[
\delta(u-v- \frac{i \hbar c}{2}) H^+_i(u-\frac{i \hbar c}{4} )
- \delta(u-v + \frac{i \hbar c}{2}) H^-_i(v-\frac{i \hbar c}{4} )
\right],\\
E_i(u_1) E_i(u_2) E_j(v)  \!\!\!\!&-&\!\!\!\! 2 \mbox{cos}(\pi \eta \hbar)
E_i(u_1) E_j(v) E_i(u_2) +  E_j(v) E_i(u_1)E_i(u_2) \nonumber \\
& &~~~~~~~~+ (u_1 \leftrightarrow u_2 ) = 0,~~~~~\mbox{for}~~A_{ij}=-1,
\label{Serr1} \\
F_i(u_1) F_i(u_2) F_j(v) \!\!\!\!&-&\!\!\!\! 2 \mbox{cos}(\pi \eta' \hbar)
F_i(u_1) F_j(v) F_i(u_2) +  F_j(v) F_i(u_1)F_i(u_2) \nonumber \\
& &~~~~~~~~+ (u_1 \leftrightarrow u_2 ) = 0,~~~~~\mbox{for}~~A_{ij}=-1,
\label{Serr2} \\
& &[c, everything] = 0 = [1, everything], \label{CurrE}
\end{eqnarray}

\noindent {\em where $u, v$ etc. are spectral parameters,
real $\hbar$, $\eta$ are two deformation parameters, $B_{ij}=A_{ij}/2$,
$A_{ij}$ are matrix elements of the Cartan matrix for the
Lie algebra $g$, and}\footnote{We assume throughout this paper
that $\eta$ and $\hbar$ are generic, i.e. $\hbar$ is not a
rational multiple of $\eta$.}

\begin{eqnarray*}
\frac{1}{\eta'} - \frac{1}{\eta} = \hbar c.
\end{eqnarray*}

\begin{rem}
For $g =sl_2$, the above current algebra reduces to the current
realization of  \Al, where the Serre-like relations (\ref{Serr1}-\ref{Serr2})
are not present.
\end{rem}

\begin{rem}
In the limit $\eta \rightarrow 0$, the current algebra \Alg
would have the same form as that of the Yangian double $DY_\hbar(g)$.
But the limiting algebra \Alo should not be considered to be
isomorphic with the Yangian double $DY_\hbar(g)$ because the
element of the algebra \Alo carries a {\em continuous} index
whilst that of the Yangian double $DY_\hbar(g)$ carries discrete one.
For $g=sl_2$, see \cite{KLP}  for more information on this point.
\end{rem}

To have a precise definition for the algebra \Alg (and {\em not} its
current realization form), we have to consider two different
cases as did in ref.\cite{KLP} for $sl_2$ case: 1) the case $c\neq 0$ and
2) the case $c=0$. In the first case one should consider the currents
$\{H^{\pm}_i(u),~E_i(u),~F_i(u)\}$ as the following Fourier transforms
of the actual elements $t_i(\lambda),~e_i(\lambda)$ and $f_i(\lambda)$
($\lambda \in \Bbb{R}$) of the algebra \Alg,

\begin{eqnarray*}
& & H^\pm_i(u) = -\frac{\hbar}{2} \int_{-\infty}^{\infty} d \lambda ~
\mbox{e}^{i \lambda u} t_i(\lambda) \mbox{e}^{\mp \lambda / 2\eta''},\\
& & E_i(u) = \int_{-\infty}^{\infty} d \lambda ~
\mbox{e}^{i \lambda u} e_i(\lambda),\\
& & F_i(u) = \int_{-\infty}^{\infty} d \lambda ~
\mbox{e}^{i \lambda u} f_i(\lambda),
\end{eqnarray*}

\noindent whereas in the second case, the currents $H^\pm_i(u)$ should
be given another expression in terms of the actual elements
$h_i(\lambda)$ of \Alg at $c=0$,

\begin{eqnarray*}
H^\pm_i(u) = \mp \hbar \int_{-\infty}^{\infty} d \lambda ~
\mbox{e}^{i \lambda u} \frac{h_i(\lambda)}{1-\mbox{e}^{\mp\lambda/
\eta}}.
\end{eqnarray*}

\noindent The difference between the cases for
$c \neq 0$ and $c=0$ can be summarized in a more compact
relationship between the algebra generators $t_i(\lambda)$
and $h_i(\lambda)$. In fact, from the two expressions of
$H^\pm_i(u)$, we can write down the following relation
at $c=0$,

\begin{eqnarray*}
h_i(\lambda) = t_i(\lambda) \mbox{sh} \left(\frac{\lambda}
{2 \eta} \right).
\end{eqnarray*}

\noindent Therefore, in the limit $c \rightarrow 0$,
$h_i(0)$ is well-defined but $t_i(0)$ tend to infinity.
On the contrary, when $c \neq 0$, $t_i(0)$ is well-defined
and $h_i(0)$ tends to zero.

Given the above Fourier transformations, one can in principle
write down the generating relations for the algebra \Alg
in terms of the continuous generators $t_i(\lambda)~(h_i(\lambda)),~
e_i(\lambda),~f_i(\lambda)$. However, such relations are
rather complicated and they are of no use in the rest of
this paper. Therefore we shall omit such relations
and consider only the {\em current realization}
(\ref{CurrB}-\ref{CurrE}) of the algebra \Alg.

Unlike the usual $q$-affine algebras and the Yangian doubles,
the algebra under consideration is {\em not} a Hopf algebra.
Recall that a Hopf algebra ${\cal A}$ is an algebra
endowed with five operations:

\begin{itemize}
\item the algebra multiplication $m: {\cal A} \times {\cal A}
\rightarrow {\cal A},~m(a \otimes b)
=ab~\mbox{for}~\forall a,b \in {\cal A}$;
\item the unit embedding $\iota: \Bbb{C} \rightarrow {\cal A}$,
$\iota(c)=c1$,~$c \in \Bbb{C}$,~$1 \in {\cal A}$ is the unit element;
\item comultiplication $\Delta: {\cal A}
\rightarrow {\cal A} \times {\cal A},~
\Delta(ab)=\Delta(a)\Delta(b)~\mbox{for}~\forall a,b
\in {\cal A}$;
\item the antipode $S: {\cal A} \rightarrow
{\cal A},~S(ab)=S(b)S(a)~\mbox{for}~\forall a,b \in {\cal A}$;
\end{itemize}

\noindent and

\begin{itemize}
\item the counit $\epsilon: {\cal A} \rightarrow \Bbb{C}$,
{}$~\epsilon(a_i)=c_i,~\mbox{for}~\forall a_i \in {\cal A}
{}~\mbox{and}~c_i \in \Bbb{C}$.
\end{itemize}

\noindent To make the algebra ${\cal A}$ into a Hopf algebra, these
structures have to obey the following axioms,

\begin{eqnarray}
& & m \circ (m \otimes \mbox{id})
= m \circ ( \mbox{id} \otimes m),\\
& & (\Delta \otimes \mbox{id}) \circ \Delta
= (\mbox{id} \otimes \Delta) \circ \Delta, \label{Del2} \\
& & (\epsilon \otimes \mbox{id}) \circ \Delta
= \mbox{id} = (\mbox{id} \otimes \epsilon) \circ \Delta,\\
& & m \circ (S \otimes \mbox{id}) \circ \Delta =
\epsilon = m \circ (\mbox{id} \otimes S) \circ \Delta.
\end{eqnarray}

\noindent For our algebra \Alg, only the first of these
axioms holds, which ensures the associativity of the algebra
multiplication. The operations $\Delta,~\epsilon,~S$ cannot
be defined on the algebra \Alg alone. However, as first discovered
in ref.\cite{KLP}, a well-defined coproduct can be defined over
the so-called ``Hopf family of algebras'' containing a finite number of
algebras of the kind \Al but with different parameters $\eta$.
However, as stated in the introduction, the case for arbitrary $g$
is much more complicated and we can only make our analysis in the current
realization. This difficulty prevented us from obtaining an
analog structure of KLP's Hopf family of algebras because the analogous
comultiplication is not known. Therefore we proceed to introduce
an alternative notion--the infinite Hopf family of algebras.
It should be remarked that no relationship is implied here between
our infinite Hopf family of algebras and the (finite) Hopf
family of algebras introduced by KLP.

\begin{defin}
Let $\{{\cal A}_n, n \in \Bbb{Z}\}$ be a family of associative
algebras with unit defined over $\Bbb{C}$. If on each ${\cal A}_n$
one can define the following operations

\begin{itemize}
\item the comultiplications $\Delta_{n}^{+}: {\cal A}_n \rightarrow
{\cal A}_n \times {\cal A}_{n + 1},~\Delta_{n}^{+}
(a_{(n)})= a_{(n)} \otimes a_{(n + 1)}$ and
$\Delta_{n}^{-}: {\cal A}_n \rightarrow
{\cal A}_{n-1} \times {\cal A}_{n},~\Delta_{n}^{-}
(a_{(n)})= a_{(n-1)} \otimes a_{(n)}$, where $a_{(n)} \in {\cal A}_n$
and $\Delta_n^{\pm}$ are algebra morphisms;
\item the counits $\epsilon_n: {\cal A}_n \rightarrow \Bbb{C}$;
\item the antipodes $S_{n}^{\pm}: {\cal A}_n \rightarrow
{\cal A}_{n \pm 1},~S_{n}^{\pm}(a_{(n)}b_{(n)}) = S_{n}^{\pm}
(b_{(n\pm 1)})S_{n}^{\pm}(a_{(n \pm 1)})$, which are algebra
anti-morphisms,
\end{itemize}

\noindent and if they satisfy the following constraints,

\begin{eqnarray}
& &(\epsilon_n \otimes \mbox{id}_{ n+1}) \Delta_{n}^{+} =
\mbox{id}_{n+1}, \label{AxB} \\
& & (\mbox{id}_{n-1} \otimes \epsilon_{n})
\Delta_{n}^{-} = \mbox{id}_{n-1},\\
& & m_{n+1} \circ (S_{n}^{+} \otimes \mbox{id}_{n+1}) \circ
\Delta_{n}^{+} = \epsilon_{n+1},\\
& & m_{n-1} \circ (\mbox{id}_{n-1} \otimes S_{n}^{-}) \circ
\Delta_{n}^{-}=\epsilon_{n-1}, \label{AxE}
\end{eqnarray}

\noindent where $m_n$ is the algebra multiplication for the
$n$-th component algebra ${\cal A}_n$, then we call the
family of algebras $\{{\cal A}_n, n\in \Bbb{Z}\}$ an infinite Hopf family of
algebras.
\end{defin}

A trivial example for  the infinite Hopf family of algebras is the
family $\{{\cal A}_n \equiv {\cal A}, n\in \Bbb{Z}\}$ in which
${\cal A}$ is a usual Hopf algebra. In this case, all our
axioms (\ref{AxB}-\ref{AxE}) hold with the comultiplications
$\Delta_{n}^{\pm}$, counits $\epsilon_n$ and the
antipodes $S_n^\pm$ being identified with those corresponding
structures for the usual Hopf algebra. Notice that
in this trivial case, we have one more axiom, eq.(\ref{Del2}),
which represents the coassociativity of the comultiplication.
For general cases no coassociativity is required. One may consider
the lost of coassociativity in our infinite Hopf family of algebras
a serious drawback compared to the (finite) Hopf family structure
of \cite{KLP}. However it will soon be clear in Preposition 2
that this structure
would bring about a great advantage in obtaining realizations of
our algebra at integer levels $k>1$.

Now we proceed to construct a nontrivial example for the
infinite Hopf family of algebras containing our algebra \Alg
as a member.

Let $\eta^{(0)} = \eta$. For all $n \in \Bbb{Z}$, let us define
$\eta^{(n)}$ recursively such that

\begin{eqnarray*}
\frac{1}{\eta^{(n+1)}} - \frac{1}{\eta^{(n)}}= \hbar c_n,
\end{eqnarray*}

\noindent where $c_n$ are a set of parameters and $c_0 \equiv c$,
the center of our algebra \Alg. Clearly, for $n=0$, we have
$\eta^{(1)}= \eta'$. The notations $\Algi{\eta^{(n)}}{c_n}$
have obvious meaning with the specification $\Algi{\eta^{(0)}}{c_0}
=$\Alg.

\begin{prop}
The family of algebras $\{{\cal A}_n \equiv \Algi{\eta^{(n)}}{c_n},
{}~n\in \Bbb{Z} \}$ form an infinite Hopf family of algebras with the
comultiplications $\Delta_n^\pm$, counits $\epsilon_n$ and antipodes
$S_n^\pm$ defined as follows,

\begin{itemize}
\item the comultiplications $\Delta_n^\pm$:
\begin{eqnarray}
\Delta_n^+ c_n \!\!\!\!&=&\!\!\!\! c_n + c_{n+1}, \label{CoB} \\
\Delta_n^- c_n \!\!\!\!&=&\!\!\!\! c_{n-1} + c_n,\\
\Delta_n^+ H_i^+(u;\eta^{(n)}) \!\!\!\!&=&\!\!\!\!
H_i^+(u+ \frac{i \hbar c_{n+1}}{4};\eta^{(n)})
\otimes H_i^+(u- \frac{i \hbar c_{n}}{4};\eta^{(n+1)}),\\
\Delta_n^- H_i^+(u;\eta^{(n)}) \!\!\!\!&=&\!\!\!\!
H_i^+(u+ \frac{i \hbar c_{n}}{4};\eta^{(n-1)})
\otimes H_i^+(u- \frac{i \hbar c_{n-1}}{4};\eta^{(n)}),\\
\Delta_n^+ H_i^-(u;\eta^{(n)}) \!\!\!\!&=&\!\!\!\!
H_i^-(u- \frac{i \hbar c_{n+1}}{4};\eta^{(n)})
\otimes H_i^-(u+ \frac{i \hbar c_{n}}{4};\eta^{(n+1)}),\\
\Delta_n^- H_i^-(u;\eta^{(n)}) \!\!\!\!&=&\!\!\!\!
H_i^-(u- \frac{i \hbar c_{n}}{4};\eta^{(n-1)})
\otimes H_i^-(u+ \frac{i \hbar c_{n-1}}{4};\eta^{(n)}),\\
\Delta_n^+ E_i(u;\eta^{(n)}) \!\!\!\!&=&\!\!\!\! E_i(u; \eta^{(n)})
\otimes 1 + H^-_i(u + \frac{i \hbar c_n}{4}; \eta^{(n)}) \otimes
E_i(u+ \frac{i \hbar c_n}{2}; \eta^{(n+1)}),\\
\Delta_n^- E_i(u;\eta^{(n)}) \!\!\!\!&=&\!\!\!\! E_i(u; \eta^{(n-1)})
\otimes 1 + H^-_i(u + \frac{i \hbar c_{n-1}}{4};
\eta^{(n-1)}) \otimes E_i(u+ \frac{i \hbar c_{n-1}}{2};
\eta^{(n)}),\\
\Delta_n^+ F_i(u;\eta^{(n)}) \!\!\!\!&=&\!\!\!\! 1 \otimes
F_i(u; \eta^{(n+1)}) + F_i(u+ \frac{i \hbar c_{n+1}}{2}; \eta^{(n)})
\otimes H^+_i(u + \frac{i \hbar c_{n+1}}{4}; \eta^{(n+1)}),\\
\Delta_n^- F_i(u;\eta^{(n)}) \!\!\!\!&=&\!\!\!\! 1 \otimes
F_i(u; \eta^{(n)}) + F_i(u+ \frac{i \hbar c_{n}}{2}; \eta^{(n-1)})
\otimes H^+_i(u + \frac{i \hbar c_{n}}{4}; \eta^{(n)}); \label{CoE}
\end{eqnarray}

\item the counits $\epsilon_n$:
\begin{eqnarray*}
\epsilon_n ( c_n )\!\!\!\!&=&\!\!\!\!0,\\
\epsilon_n ( 1_n ) \!\!\!\!&=&\!\!\!\! 1,\\
\epsilon_n ( H^\pm_i(u; \eta^{(n)}))  \!\!\!\!&=&\!\!\!\! 1,\\
\epsilon_n ( E_i(u; \eta^{(n)})) \!\!\!\! &=& \!\!\!\! 0,\\
\epsilon_n ( F_i(u; \eta^{(n)})) \!\!\!\! &=& \!\!\!\! 0;
\end{eqnarray*}

\item the antipodes $S_n^\pm$:

\begin{eqnarray*}
S_n^\pm (c_n) \!\!\!\!&=&\!\!\!\! - c_{n \pm 1},\\
S_n^\pm (H^\pm_i(u; \eta^{(n)}))\!\!\!\!& =&\!\!\!\!
[ H^\pm_i(u; \eta^{(n \pm 1)} )]^{-1},\\
S^\pm_n (E_i(u;\eta^{(n)})) \!\!\!\!&=&\!\!\!\!
- H^-_i(u - \frac{i \hbar c_{n \pm 1}}{4}; \eta^{(n \pm 1)})^{-1}
E_i(u- \frac{i \hbar c_{n \pm 1}}{2}; \eta^{(n \pm 1)}),\\
S^\pm_n(F_i(u;\eta^{(n)})) \!\!\!\!&=&\!\!\!\!
- F_i(u - \frac{i \hbar c_{n \pm 1}}{2}; \eta^{(n \pm 1)})
H^+_i(u- \frac{i \hbar c_{n \pm 1}}{4}; \eta^{(n \pm 1)})^{-1}.
\end{eqnarray*}
\end{itemize}

\noindent where the second arguments in the current operators
(the $\eta$'s) indicate to which algebra the currents belong.
\end{prop}

The proof for this proposition is straightforward. Notice that
in this example the comultiplications $\Delta^\pm_n$ are not all
independent. A simple observation would show that

\begin{eqnarray*}
\Delta^-_n{\cal A}_n = \Delta^+_{n-1} {\cal A}_{n-1}.
\end{eqnarray*}

\noindent Two more remarks are in due course.

\begin{rem}
In the case of $c_n=0$ for all $n\in \Bbb{Z}$, the infinite Hopf family of algebras
become trivial again because there are no differences between the algebras
$\Algi{\eta^{(n)}}{0}$ and $\Algi{\eta^{(m)}}{0}$ for any pair of
$n,m \in \Bbb{Z}$.
\end{rem}

\begin{rem} Under the cases of remarks 2 and 3, the above structures
for the infinite Hopf family of algebras reduce to the original Hopf algebra
structure. In particular, under the case of remark 2, the
comultiplications would have the same form with the so-called Drinfeld
comultiplication for the Yangian double.
\end{rem}

The comultiplications introduced above are useful not only
in clarifying the structure of the infinite Hopf family of algebras
but also in the representation theory of the
representative algebra \Alg. Before going into detailed
structure of representations, we state the following
proposition, which can be directly verified.

\begin{prop} \label{zxcc}
The comultiplication $\Delta^+_n$
defined in eqs.(\ref{CoB}-\ref{CoE}) induce algebra homomorphism from
$\Algi{\eta^{(n)}}{c_n} \otimes
\Algi{\eta^{(n+1)}}{c_{n+1}}$ to
$\Algi{\eta^{(n)}}{c_n + c_{n+1}}$, $\Delta^-_n$ induce
homomorphism from \\
$\Algi{\eta^{(n-1)}}{c_{n-1}} \otimes
\Algi{\eta^{(n)}}{c_{n}}$ to
$\Algi{\eta^{(n-1)}}{c_{n-1} + c_{n}}$.
\end{prop}

\noindent Actually, the above proposition states that
the images of the generating currents $E_i(u; \eta^{(n)}),
{}~F_i(u; \eta^{(n)})$ and $H^\pm_i(u; \eta^{(n)})$
of $\Algi{\eta^{(n)}}{c_n}$ under $\Delta^\pm_n$ satisfy the
defining relations for $\Algi{\eta^{(n)}}{c_n + c_{n+1}}$
and $\Algi{\eta^{(n-1)}}{c_{n-1} + c_{n}}$ respectively.
This result is quite astonishing at on hand,
and will be quite useful for constructing a higher level
realization out of level 1 representations on the other.
Therefore we proceed to consider the level 1 representation of our algebra.

\section{Representation theory}

\subsection{Free boson realization of \Alg at level $c=1$}

First we would like to consider the free boson realization of the
generating relations (\ref{CurrB}-\ref{CurrE}) for the algebra \Alg.
For this we introduce the set of deformed free bosons $a_i(\lambda)$
with continuous parameter $\lambda \neq 0$ and discrete
$i=1,...,r$, which constitute the following deformed
Heisenberg algebra ${\cal H}(\eta)$:

\begin{equation}
[ a_i(\lambda), a_j(\mu) ] = \frac{4}{\lambda} \mbox{sh}
\frac{\hbar\lambda}{2} \mbox{sh}(\hbar B_{ij} \lambda)
\frac{\mbox{sh}\frac{\lambda}{2\eta}}{\mbox{sh}
\frac{\lambda}{2\eta'}} \delta(\lambda + \mu). \label{Heis}
\end{equation}

\noindent We also use the notations $a'_i(\lambda) =
a_i(\lambda) \frac{\mbox{sh}\frac{\lambda}{2\eta'}}{\mbox{sh}
\frac{\lambda}{2\eta}}$, which satisfy the relations

\begin{eqnarray*}
[ a_i'(\lambda), a_j'(\mu) ] = \frac{4}{\lambda} \mbox{sh}
\frac{\hbar\lambda}{2} \mbox{sh}(\hbar B_{ij} \lambda)
\frac{\mbox{sh}\frac{\lambda}{2\eta'}}{\mbox{sh}
\frac{\lambda}{2\eta}} \delta(\lambda + \mu).
\end{eqnarray*}

The normal ordering for the exponential expressions of the
above free bosons are defined in the following way \cite{KLP},

\begin{eqnarray}
& & :\mbox{exp} \int_{-\infty}^{\infty} d\lambda~g_1(\lambda)a_i(\lambda):
:\mbox{exp} \int_{-\infty}^{\infty} d\mu~g_2(\mu)a_j(\mu):  \nonumber \\
& &~~~~= \mbox{exp}\left( \int_C\frac{d\lambda \mbox{ln}(-\lambda)}{2\pi i}
\alpha_{ij}(\lambda) g_1(\lambda) g_2(- \lambda) \right)
:\mbox{exp} \left( \int_{-\infty}^{\infty} d\lambda~
g_1(\lambda)a_i(\lambda) + \int_{-\infty}^{\infty} d\mu~
g_2(\mu) a_j(\mu) \right):, \nonumber\\
\label{norm}
\end{eqnarray}

\noindent where $\alpha_{ij}(\lambda)$ is a function given by

\begin{equation}
[ a_i(\lambda), a_j(\mu) ] = \alpha_{ij}(\lambda) \delta(\lambda + \mu),
\label{AbbHeis}
\end{equation}

\noindent $C$ is a contour on the complex $\lambda$-plane
depicted in Figure 1. Moreover, we introduce the following zero mode operators,

\begin{eqnarray*}
[P_i, Q_j] =B_{ij}.
\end{eqnarray*}

\begin{prop}
The following bosonic expressions realize the generating relations
(\ref{CurrB}-\ref{CurrE}) of the algebra \Alg with $c=1$,
\end{prop}

\begin{eqnarray}
E_j(u) &=& \mbox{e}^{\gamma}~\mbox{exp}(2\pi i Q_j)~ \mbox{exp}(P_j)
:\mbox{exp}\left( \frac{1}{2} \phi'_j(u) \right):, \label{BosB} \\
F_j(u) &=& \mbox{e}^{\gamma}~\mbox{exp}(-2\pi i Q_j)~ \mbox{exp}(-P_j)
:\mbox{exp}\left(- \frac{1}{2} \phi_j(u) \right):,\\
H^\pm_j(u)&=&\mbox{e}^{- 2\gamma}~
:E_i(u \pm \frac{i \hbar}{4}) F_j(u \mp \frac{i\hbar}{4}):\\
&=& :\mbox{exp} \left( \mp \int_{-\infty}^{\infty}
d\lambda~\mbox{e}^{i\lambda u} \frac{\mbox{e}^{\mp \hbar \lambda /4}}{
1-\mbox{e}^{\pm \lambda/\eta}} a_j(\lambda)\right):,
\end{eqnarray}

\noindent where

\begin{eqnarray}
& & \phi_j(u)=\int_{-\infty}^{+\infty}d \lambda~\mbox{e}^{i\lambda
u} \frac{a_j(\lambda)}{\mbox{sh}\frac{\hbar\lambda}{2}},\\
& & \phi'_j(u)=\int_{-\infty}^{+\infty}d \lambda~\mbox{e}^{i\lambda
u} \frac{a'_j(\lambda)}{\mbox{sh}\frac{\hbar\lambda}{2}}. \label{BosE}
\end{eqnarray}

The proof of this proposition is also by straightforward but tedious
calculations. The normal ordering rule (\ref{norm}) and the following formula
which can be found in ref.\cite{KLP} are very useful for the calculations,

\begin{eqnarray*}
\int_C \frac{d\lambda~\mbox{ln}(-\lambda)}{2\pi i \lambda}
\frac{\mbox{e}^{-x\lambda}}{1-\mbox{e}^{-\lambda/ \eta}} = \mbox{ln}
\Gamma(\eta x) + (\eta x -\frac{1}{2}) (\gamma-\mbox{ln} \eta)
-\frac{1}{2} \mbox{ln} (2\pi),
\end{eqnarray*}

\noindent where $\Gamma(x)$ is the usual Gamma function which satisfy the
following formula,

\begin{eqnarray*}
\Gamma(x)\Gamma(1-x) = \frac{\pi}{\mbox{sin}\pi x}.
\end{eqnarray*}

\vspace{0.5cm}
\epsfbox{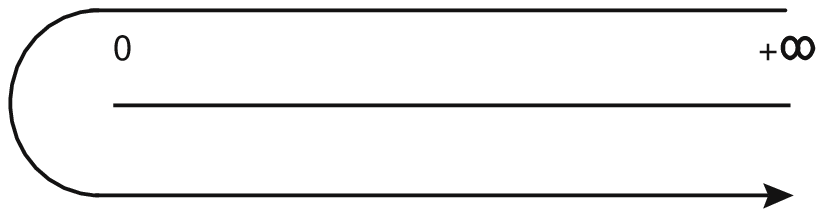}
\begin{center}
{Figure 1: The integration contour $C$}
\end{center}

It is interesting to mention that the bosonization formulas for
the currents $E_i(u)$ and $F_i(u)$ are quite similar to that
of the screening currents of the quantum $(\hbar,\xi)$-deformed
$W$-algebras \cite{YWL}.

\subsection{Representations at other integer levels}

The bosonic expressions (\ref{BosB}-\ref{BosE}) only give a bosonic
realization of \Alg at level $c=1$. However, as mentioned in the last section,
it is possible to obtain realizations of \Alg at other integer levels
using the knowledge gathered so far. The key point is to make use of
Proposition 2 repeatedly, first in the case of $c=c_{0}=c_{1}=1$ (which lead to
a realization at level $c=2$),
then in the case of $c=c_{0}=2, c_{1}=1$, and so on.

We give the following proposition

\begin{prop}
The level $c=k~(k \in \Bbb{Z}_+)$ bosonic realization for the algebra \Alg
can be obtained using $k$ copies of the Heisenberg algebra (\ref{Heis})
$ \{ {\cal H}(\eta^{(l)}), l=0,~1,~...,~k-1 \}$ (each of which
realizes the level 1 representation for the algebras
$\Algi{\eta^{(l)}}{}$ with $l=0,~1,~...,~k-1$) and the repeated use of
Proposition 2 (the comultiplication $\Delta^{+}_0$).
\end{prop}

Actually, the above proposition provides a way to understand the meaning
of the infinite Hopf family of algebras--instead of getting higher level
representations of any distinguished member of this family, one can
study the level 1 representations for several members simultaneously.

To obtain bosonic realizations of negative integer level,
one may use the antipodes $S^{\pm}_n$. However such realizations
are of less interests to us.

\subsection{The structure of Fock spaces}

The bosonic realizations of the algebra \Alg can be viewed
as representations on the Fock space of the bosonic Heisenberg
algebras. Therefore, for completeness, we have to pay some words
on the structure of Fock spaces.

First we specify the Fock space for the level 1 representation.
Consider the abbreviated form (\ref{AbbHeis}) of the bosonic
Heisenberg algebra ${\cal H}(\eta)$. The structure functions
$\alpha_{ij}(\lambda)$ have the properties

 \begin{eqnarray}
 & & \alpha_{ij}(\lambda) = - \alpha_{ij}(-\lambda), \nonumber \\
 & & \alpha_{ij}(\lambda) = \alpha_{ji}(\lambda). \label{sym}
 \end{eqnarray}

\noindent Let $|vac\rangle$ be a right ``vacuum state''.
The right Fock space ${\cal F}(\eta)$ is generated from
$|vac\rangle$ as follows,

\begin{eqnarray*}
\int_{-\infty}^0
d\lambda_n~f_n(\lambda_n) a_{i_n}(\lambda_n)
... \int_{-\infty}^0 d\lambda_1~f_1(\lambda_1) a_{i_1}(\lambda_1)
~| vac \rangle,~~~
i_l=1,~2,~...,~r,~~~l=1,~2,~...,~n,
\end{eqnarray*}

\noindent where $f_l(\lambda)$ are functions which
are analytic in a neighborhood of $\Bbb{R}_+$ except
$\lambda=0$, where a simple pole may appear. For each concrete
$\alpha_{ij}(\lambda)$, proper asymptotic behaviors for
$f_l(\lambda)$ as $\lambda \rightarrow +\infty$ are
required. However, we do not specify them in detail (for the special
case of $g=sl_2$, such asymptotics were given explicitly in
ref.\cite{KLP}).

Similarly, let $\langle vac|$ be a left ``vacuum state''.
The left Fock space ${\cal F}^\ast(\eta)$ is generated from
as follows,

\begin{eqnarray*}
\langle vac|~\int^{+\infty}_0
d\lambda_1~g_1(\lambda_1) a_{i_1}(\lambda_1)
... \int^{+\infty}_0 d\lambda_n~g_n(\lambda_n) a_{i_n}(\lambda_n),
{}~~~i_l=1,~2,~...,~r,~~~l=1,~2,~...,~n,
\end{eqnarray*}

\noindent where $g_l(\lambda)$ are functions which
are analytic in a neighborhood of $\Bbb{R}_-$ except
$\lambda=0$, where a simple pole may appear. As in the case of
$f_l(\lambda)$, proper asymptotic behaviors for
$g_l(\lambda)$ as $\lambda \rightarrow -\infty$ are also
required.

Like in the case of $g=sl_2$, the pairing $(~,~):{\cal F}^\ast(\eta) \times
{\cal F}(\eta) \rightarrow \Bbb{C}$ between the left and right
Fock spaces can be uniquely defined by the following prescriptions,

\begin{itemize}
\item $(\langle vac|,|vac \rangle)=1,$
\item $ \displaystyle
(\langle vac|~\int^{+\infty}_0
d\lambda~g(\lambda) a_{i}(\lambda),~\int_{-\infty}^0
d\mu~f(\mu) a_{j}(\mu)~|vac \rangle)
=\int_C \frac{d\lambda\mbox{ln}(-\lambda)}{2\pi i}g(\lambda)
f(-\lambda) \alpha_{ij}(\lambda),$
\item the Wick theorem.
\end{itemize}

Now let the vacuum states $|vac \rangle$ and $\langle vac|$ be such that

\begin{eqnarray*}
& &a_i(\lambda)~|vac \rangle =0,~~~\lambda >0,~~~~~~~ P_i~|vac \rangle =0,\\
& &\langle vac|~a_i(\lambda) =0,~~~\lambda<0,~~~~~~~ \langle vac|~Q_i =0.
\end{eqnarray*}

\noindent Let $f(\lambda)$ be analytic in some neighborhood of the real
$\lambda$-line, satisfying proper analytic behaviors as
$\lambda \rightarrow \pm \infty$, and may have simple poles at $\lambda=0$.
Then the action of the expressions like

\begin{eqnarray*}
F=:\mbox{exp}\left( \int_{-\infty}^{+\infty} d\lambda~f(\lambda)
a_i(\lambda) \right):
\end{eqnarray*}

\noindent on ${\cal F}(\eta)$ and ${\cal F}^\ast(\eta)$ are
given respectively by the decompositions $F=F_-F_+$ and
$F=\tilde{F}_-\tilde{F}_+$, where

\begin{eqnarray*}
& & F_-=\mbox{exp}\left( \int_{-\infty}^{0} d\lambda~f(\lambda)
a_i(\lambda) \right),\\
& & F_+=\lim_{\epsilon \rightarrow 0^+}
\mbox{e}^{\epsilon\mbox{ln} \epsilon f(\epsilon) a_i(\epsilon)}
\mbox{exp}\left( \int_{\epsilon}^{+\infty} d\lambda~f(\lambda)
a_i(\lambda) \right),\\
& & \tilde{F}_-=\lim_{\epsilon \rightarrow 0^+}
\mbox{e}^{\epsilon\mbox{ln} \epsilon f(-\epsilon) a_i(-\epsilon)}
\mbox{exp}\left( \int^{-\epsilon}_{-\infty} d\lambda~f(\lambda)
a_i(\lambda) \right),\\
& & \tilde{F}_+=\mbox{exp}\left( \int^{+\infty}_{0} d\lambda~f(\lambda)
a_i(\lambda) \right).
\end{eqnarray*}

\noindent Moreover, these two actions are adjoint to each other,
and the product of normal ordered operators like $F$ satisfy our
normal ordering rule (\ref{norm}). This complete the description of Fock
spaces at level 1.

The Fock spaces for level $k$ bosonic representation of \Alg
is nothing but the direct product of $k$ copies of the level 1
Fock spaces, namely, ${\cal F}^{(k)}(\eta^{(0)},~...,~\eta^{(k-1)})
={\cal F}(\eta^{(0)}) \otimes ... \otimes {\cal F}(\eta^{(k-1)})$.
The left Fock space for level $k$ bosonic
representation has a similar structure, ${\cal F}^{\ast(k)}
(\eta^{(0)},~...,~\eta^{(k-1)})={\cal F}^\ast(\eta^{(0)})
\otimes ... \otimes {\cal F}^\ast(\eta^{(k-1)})$.

\subsection{The case of $c=0$: evaluation representation}

As mentioned earlier, the structure of the algebra \Alg changes
drastically from $c\neq 0$ to $c=0$. This change is not only reflected
in the different asymptotic behaviors for the generating currents,
but also in the trivialization of the structure of the infinite Hopf family
(see Remark 3), and it also affects the representation theory
at $c=0$.

Just like the usual affine Lie algebras and the  affine Hopf algebras,
among the class of level 0 representations for the algebra \Alg,
there is a special subclass which is finite dimensional.
We adopt the terminology from the representation theory
of affine and affine Hopf algebras and call the finite dimensional
level 0 representations of \Alg the evaluation representations.

Recall that there is no differences between the algebras
$\Algi{\eta^{(m)}}{0}$ and $\Algi{\eta^{(n)}}{0}$ for different $n$
and $m$. Recall also that the evaluation representations
for the usual affine Hopf algebras are best written in terms of
``half currents'' rather than the total currents which we have been using
for \Alg so far. Therefore it seem that the first step to give an
evaluation representation for the algebra \Alg is to split the total
currents $E_i(u)$ and $F_i(u)$ into half currents.
This task can be fulfilled in a completely analogous way as in the $sl_2$
case.

We define (for generic $c$) the half currents $e^\pm_i(u)$
and $f^\pm_i(u)$ as follows,

\begin{eqnarray*}
& &e^\pm_i(u) = \pi \eta \int_{C_1} \frac{dv}{2\pi i}
\frac{E_i(v)}{\mbox{sh} \pi\eta(u-v\pm i\hbar c/4)},\\
& &f^\pm_i(u) = \pi \eta' \int_{C_2} \frac{dv}{2\pi i}
\frac{F_i(v)}{\mbox{sh} \pi\eta'(u-v\mp i\hbar c/4)},
\end{eqnarray*}

\noindent where the contours $C_1$ and $C_2$ run from $-\infty$ to $\infty$,
with the points $u+i\hbar c/4 - ik/\eta(k \geq 0)$ above $C_1$,
$u+i\hbar c/4 + ik/\eta~(k \geq 0)$ below $C_1$,
$u+i\hbar c/4 + ik/\eta'~(k \geq 0)$ above $C_2$,
$u-i\hbar c/4 + ik/\eta'~(k \geq 0)$ below $C_2$.

The remarkable point for these half currents is that they
satisfy the following Ding-Frenkel like relations,

\begin{eqnarray*}
& & e^+_i(u-\frac{i\hbar c}{4})-e^-_i(u+ \frac{i\hbar c}{4}) =E_i(u),\\
& & f^+_i(u+\frac{i\hbar c}{4})-f^-_i(u- \frac{i\hbar c}{4}) =F_i(u),
\end{eqnarray*}

\noindent however these relations should be understood in some
proper analytic continuation sense in contrast to the direct
decompositions of formal power series \cite{iohara,KT}. To be explicit, we give
the domains of analyticity for the half currents:

\begin{eqnarray*}
& & e^+_i(u),~f^+_i(u),~H^+_i(u):~~\mbox{analytic in}~~\Pi_+=
\left\{ -\frac{1}{\eta}-\frac{\hbar c}{4} < \mbox{Im} u
< - \frac{\hbar c}{4} \right\},\\
& & e^-_i(u),~f^-_i(u),~H^-_i(u):~~\mbox{analytic in}~~\Pi_-=
\left\{ \frac{\hbar c}{4} < \mbox{Im} u
< \frac{\hbar c}{4}+ \frac{1}{\eta} \right\}.
\end{eqnarray*}

\noindent Moreover, we have

\begin{eqnarray*}
e^-_i(u)=-e^+_i(u-i/\eta''),~~~f^-_i(u)=-f^+_i(u-i/\eta'')
\end{eqnarray*}

\noindent and

\begin{eqnarray*}
H^-_i(u)=H^+_i(u-i/ \eta''),
\end{eqnarray*}

\noindent where $u \in \Pi_-$.

The following proposition gives a simplest evaluation representation
for \Alg with $g=sl_{r+1}$.

\begin{prop}
Let $V$ be an $(r+1)$-dimensional vector space with orthogonal basis
$\{v_0,~v_1,~...,~v_r\}$. The $(r+1)$-dimensional evaluation representation
of \Alg with $g=sl_{r+1}$ on $V_z(\eta)=V \otimes$ \\
$\Bbb{C}$ $[[ \mbox{e}^{\pi\eta z} ]]$ is given by
the following actions $(u \in \Pi_+)$,
\end{prop}

\begin{eqnarray*}
e^+_l(u)v_{j,z} &=& \delta_{lj} \frac{-\mbox{sh} i\pi \eta \hbar}{
\mbox{sh} \pi \eta (u-z-\frac{r-l}{2}i\hbar)} v_{j-1,z},\\
f^+_l(u)v_{j-1,z} &=& \delta_{lj} \frac{-\mbox{sh} i\pi \eta \hbar}{
\mbox{sh} \pi \eta (u-z-\frac{r-l}{2}i\hbar)} v_{j,z},\\
H^+_l(u)v_{j,z} &=& \delta_{lj}
\frac{\mbox{sh} \pi \eta(u-z-\frac{r-l-2}{2}i \hbar)}{
\mbox{sh} \pi \eta (u-z-\frac{r-l}{2}i\hbar)} v_{j,z}\\
&+& \delta_{l-1,j}
\frac{\mbox{sh} \pi \eta(u-z-\frac{r-l+2}{2}i \hbar)}{
\mbox{sh} \pi \eta (u-z-\frac{r-l}{2}i\hbar)} v_{j,z}\\
&+&(1-\delta_{lj}-\delta_{l-1,j})v_{j,z}.
\end{eqnarray*}

\noindent The relations for the ``negative'' half currents
are given by the same formulas but with $u \in \Pi_-$.

Notice that for $r=1$, the above evaluation
representation reduces to the one presented in \cite{KLP} for
\Al; for $\eta \rightarrow 0$, it reduces to the
$(r+1)$-dimensional evaluation
representation for $DY_\hbar(sl_{r+1})$ \cite{iohara}.

\subsection{The intertwining relations and vertex operators}

One of the important ingredients in the representation theory of
affine algebras is the intertwining operators which intertwine
the infinite-dimensional representation and its tensor product
with evaluation representation. For the infinite Hopf family of algebras,
we can define analogous objects, also called intertwining operators,
however acting on the space of tensor product of the infinite dimensional
representation of one member of the family and the evaluation representation
of the subsequent member of the same family, or on the space of tensor product
of the evaluation representation and some infinite dimensional representation
of a fixed member of the family.

Taking as the infinite dimensional representation the level 1 
bosonic representation,
as the evaluation representation the $(r+1)$-dimensional representation
obtained above for $g=sl_{r+1}$, we now proceed to give the definition of a particular
set of intertwining operators.

\begin{defin}
The intertwining operators (vertex operators) (here $\eta'=1/(\hbar+\frac{1}{\eta})$)

\begin{eqnarray*}
\begin{array}{ll}
$$\Phi(z): {\cal F}(\eta) \rightarrow {\cal F}(\eta) \otimes V_z(\eta'),$$ &
$$\Phi^\ast(z): {\cal F}(\eta) \otimes V_z(\eta') \rightarrow {\cal F}(\eta),$$ \cr
$$\Psi^\ast(z): V_z \otimes {\cal F}(\eta) \rightarrow {\cal F}(\eta),$$ &
$$\Psi(z): {\cal F}(\eta)\rightarrow V_z \otimes {\cal F}(\eta)$$
\end{array}
\end{eqnarray*}

\noindent are those commute with the action of \Alg
\footnote{In ref.\cite{KLP}, a twisted version of the vertex operators was
defined so that the commutation relations of the twisted vertex
operators yield the two-body $S$-matrix for Sine-Gordon model. In
our case, we do not have such motivations to define twisted
vertex operators. Moreover, remember that the comultiplications
of ref.\cite{KLP} is different from the one we are using.},

\begin{eqnarray*}
\begin{array}{ll}
$$\Phi(z)x=\Delta(x)\Phi(z)$$ &
$$\Phi^\ast(z) \Delta(x)=x\Phi^\ast(z)$$ \cr
$$\Psi^\ast(z) \Delta(x)=x\Psi^\ast(z)$$ &
$$\Psi(z)x=\Delta(x)\Psi(z),$$
\end{array}
\end{eqnarray*}

\noindent where $x \in$\Alg.
\end{defin}

The components of these vertex operators are defined as follows,

\begin{eqnarray*}
& & \Phi(z)v = \sum_{j=0}^r \Phi_j(z) v \otimes v_j,\\
& & \Phi^\ast(z) (v \otimes v_j) = \Phi^\ast_j(z) v,\\
& & \Psi^\ast(z) (v_j \otimes v) = \Psi^\ast_j(z) v,\\
& & \Psi(z)v = \sum_{j=0}^r v_j \otimes \Psi_j(z) v ,
\end{eqnarray*}

\noindent where $v \in {\cal F}(\eta)$ and $v_j \in V$.

Using the explicit form of the evaluation representation
given in the last subsection and the comultiplication
formulas (\ref{CoB}-\ref{CoE}), we are ready to obtain the following
intertwining relations (the commutation relations between
vertex operators and the generating currents for \Alg),

\begin{itemize}
\item Relations for $\Phi(z)$:

\begin{eqnarray*}
& & \left\{
\begin{array}{l}
$$\displaystyle
\Phi_j(z)H^+_j(u)=\frac{\mbox{sh}\pi \eta'(u-z-
\frac{r-j-2}{2}i\hbar -\frac{3}{4}i\hbar)}{
\mbox{sh}\pi \eta'(u-z-
\frac{r-j}{2}i\hbar -\frac{3}{4}i\hbar)}H^+_j(u) \Phi_j(z),$$\cr
$$\displaystyle
\Phi_{j-1}(z)H^+_j(u)=\frac{\mbox{sh}\pi \eta'(u-z-
\frac{r-j+2}{2}i\hbar -\frac{3}{4}i\hbar)}{
\mbox{sh}\pi \eta'(u-z-
\frac{r-j}{2}i\hbar -\frac{3}{4}i\hbar)}H^+_j(u) \Phi_{j-1}(z),$$\cr
$$\displaystyle
\Phi_{l}(z)H^+_j(u)=H^+_j(u) \Phi_{l}(z),$$~~{\rm otherwise}
\end{array}
\right.\\
& & \displaystyle
\left\{
\begin{array}{l}
$$\displaystyle
\Phi_j(z)H^-_j(u)=\frac{\mbox{sh}\pi \eta'(u-z-
\frac{r-j-2}{2}i\hbar -\frac{1}{4}i\hbar)}{
\mbox{sh}\pi \eta'(u-z-
\frac{r-j}{2}i\hbar -\frac{1}{4}i\hbar)}H^-_j(u) \Phi_j(z),$$\cr
$$\displaystyle
\Phi_{j-1}(z)H^-_j(u)=\frac{\mbox{sh}\pi \eta'(u-z-
\frac{r-j+2}{2}i\hbar -\frac{1}{4}i\hbar)}{
\mbox{sh}\pi \eta'(u-z-
\frac{r-j}{2}i\hbar -\frac{1}{4}i\hbar)}H^-_j(u) \Phi_{j-1}(z),$$\cr
$$\Phi_{l}(z)H^-_j(u)=H^-_j(u) \Phi_{l}(z),$$~~{\rm otherwise}
\end{array}
\right.\\
& & \displaystyle
[ \Phi_j(z), E_l(u) ] =\frac{\mbox{sh}i\pi \eta'\hbar}{\pi \eta'}
\delta_{j,l-1} \delta(u-z-\frac{r-l}{2}i\hbar)
H^-_l(u+\frac{i\hbar}{4}) \Phi_l(z),\\
& &\left\{
\begin{array}{l}
$$\displaystyle
\Phi_j(z)F_j(u)=\frac{\mbox{sh}\pi \eta'(u-z-
\frac{r-j-2}{2}i\hbar -\frac{i\hbar}{2})}{
\mbox{sh}\pi \eta'(u-z-
\frac{r-j}{2}i\hbar -\frac{i\hbar}{2})}F_j(u) \Phi_j(z)$$\cr
$$\displaystyle
~~~~~~~~~~~~~~~+ \frac{\mbox{sh}i\pi \eta' \hbar}{\pi \eta'}
\delta(u-z-\frac{r-l}{2}i\hbar - \frac{i\hbar}{2}) \Phi_{j-1}(z),$$\cr
$$\displaystyle
\Phi_{j-1}(z)F_j(u)=\frac{\mbox{sh}\pi \eta'(u-z-
\frac{r-j+2}{2}i\hbar - \frac{i\hbar}{2})}{
\mbox{sh}\pi \eta'(u-z-
\frac{r-j}{2}i\hbar -\frac{i\hbar}{2})}F_j(u) \Phi_{j-1}(z),$$\cr
$$\Phi_{l}(z)F_j(u)=F_j(u) \Phi_{l}(z),$$~~{\rm otherwise;}
\end{array}
\right.
\end{eqnarray*}

\item Relations for $\Phi^\ast(z)$:

\begin{eqnarray*}
& & \left\{
\begin{array}{l}
$$\displaystyle
H^+_j(u)\Phi_j^\ast(z)=\frac{\mbox{sh}\pi \eta'(u-z-
\frac{r-j-2}{2}i\hbar -\frac{3}{4}i\hbar)}{
\mbox{sh}\pi \eta'(u-z-
\frac{r-j}{2}i\hbar -\frac{3}{4}i\hbar)}\Phi_j^\ast(z)H^+_j(u),$$\cr
$$\displaystyle
H^+_j(u)\Phi_{j-1}^\ast(z)=\frac{\mbox{sh}\pi \eta'(u-z-
\frac{r-j+2}{2}i\hbar -\frac{3}{4}i\hbar)}{
\mbox{sh}\pi \eta'(u-z-
\frac{r-j}{2}i\hbar -\frac{3}{4}i\hbar)}\Phi_{j-1}^\ast(z) H^+_j(u),$$\cr
$$\displaystyle
H^+_j(u)\Phi_{l}^\ast(z)=\Phi_{l}^\ast(z)H^+_j(u),$$~~{\rm otherwise}
\end{array}
\right.\\
& & \displaystyle
\left\{
\begin{array}{l}
$$\displaystyle
H^-_j(u)\Phi_j^\ast(z)=\frac{\mbox{sh}\pi \eta'(u-z-
\frac{r-j-2}{2}i\hbar -\frac{1}{4}i\hbar)}{
\mbox{sh}\pi \eta'(u-z-
\frac{r-j}{2}i\hbar -\frac{1}{4}i\hbar)}\Phi_j^\ast(z) H^-_j(u),$$\cr
$$\displaystyle
H^-_j(u) \Phi_{j-1}^\ast(z)=\frac{\mbox{sh}\pi \eta'(u-z-
\frac{r-j+2}{2}i\hbar -\frac{1}{4}i\hbar)}{
\mbox{sh}\pi \eta'(u-z-
\frac{r-j}{2}i\hbar -\frac{1}{4}i\hbar)}\Phi_{j-1}^\ast(z) H^-_j(u),$$\cr
$$H^-_j(u) \Phi_{l}^\ast(z)=\Phi_{l}^\ast(z) H^-_j(u),$$~~{\rm otherwise}
\end{array}
\right.\\
& & \displaystyle
{}[ E_l(u), \Phi^\ast_j(z) ] =
\frac{\mbox{sh}i\pi \eta' \hbar}{\pi \eta'}
\delta_{lj} \delta(u-z-\frac{r-l}{2}i\hbar)
\Phi^\ast_{l-1}(z)H^-_l(u+\frac{i\hbar}{4}),\\
& &
\left\{
\begin{array}{l}
$$\displaystyle
F_j(u) \Phi^\ast_j(z)=\frac{\mbox{sh}\pi \eta'(u-z-
\frac{r-j-2}{2}i\hbar -\frac{i\hbar}{2})}{
\mbox{sh}\pi \eta'(u-z-
\frac{r-j}{2}i\hbar -\frac{i\hbar}{2})} \Phi^\ast_j(z) F_j(u),$$\cr
$$\displaystyle
F_j(u) \Phi^\ast_{j-1}(z)=\frac{\mbox{sh}\pi \eta'(u-z-
\frac{r-j+2}{2}i\hbar - \frac{i\hbar}{2})}{
\mbox{sh}\pi \eta'(u-z-
\frac{r-j}{2}i\hbar -\frac{i\hbar}{2})} \Phi^\ast_{j-1}(z)F_j(u)$$\cr
$$\displaystyle
~~~~~~~~~~~~~~~+ \frac{\mbox{sh}i\pi \eta' \hbar}{\pi \eta'}
\delta(u-z-\frac{r-l}{2}i\hbar -\frac{i\hbar}{2})
\Phi^\ast_{j}(z),$$\cr
$$F_j(u) \Phi^\ast_{l}(z))=\Phi^\ast_{l}(z) F_j(u),$$~~{\rm otherwise;}
\end{array}
\right.
\end{eqnarray*}

\item Relations for $\Psi^\ast(z)$:

\begin{eqnarray*}
& & \left\{
\begin{array}{l}
$$\displaystyle
H^+_j(u)\Psi_j^\ast(z)=\frac{\mbox{sh}\pi \eta(u-z-
\frac{r-j-2}{2}i\hbar -\frac{1}{4}i\hbar)}{
\mbox{sh}\pi \eta(u-z-
\frac{r-j}{2}i\hbar -\frac{1}{4}i\hbar)}\Psi_j^\ast(z)H^+_j(u),$$\cr
$$\displaystyle
H^+_j(u)\Psi_{j-1}^\ast(z)=\frac{\mbox{sh}\pi \eta(u-z-
\frac{r-j+2}{2}i\hbar -\frac{1}{4}i\hbar)}{
\mbox{sh}\pi \eta(u-z-
\frac{r-j}{2}i\hbar -\frac{1}{4}i\hbar)}\Psi_{j-1}^\ast(z) H^+_j(u),$$\cr
$$\displaystyle
H^+_j(u)\Psi_{l}^\ast(z)=\Psi_{l}^\ast(z)H^+_j(u),$$~~{\rm otherwise}
\end{array}
\right.\\
& & \displaystyle
\left\{
\begin{array}{l}
$$\displaystyle
H^-_j(u)\Psi_j^\ast(z)=\frac{\mbox{sh}\pi \eta(u-z-
\frac{r-j-2}{2}i\hbar -\frac{3}{4}i\hbar)}{
\mbox{sh}\pi \eta(u-z-
\frac{r-j}{2}i\hbar -\frac{3}{4}i\hbar)}\Psi_j^\ast(z) H^-_j(u),$$\cr
$$\displaystyle
H^-_j(u) \Psi_{j-1}^\ast(z)=\frac{\mbox{sh}\pi \eta(u-z-
\frac{r-j+2}{2}i\hbar -\frac{3}{4}i\hbar)}{
\mbox{sh}\pi \eta(u-z-
\frac{r-j}{2}i\hbar -\frac{3}{4}i\hbar)}\Psi_{j-1}^\ast(z) H^-_j(u),$$\cr
$$H^-_j(u) \Psi_{l}^\ast(z)=\Psi_{l}^\ast(z) H^-_j(u),$$~~{\rm otherwise}
\end{array}
\right.\\
& &\displaystyle
\left\{
\begin{array}{l}
$$\displaystyle
E_j(u) \Psi^\ast_j(z)=\frac{\mbox{sh}\pi \eta(u-z-
\frac{r-j-2}{2}i\hbar -\frac{i\hbar}{2})}{
\mbox{sh}\pi \eta(u-z-
\frac{r-j}{2}i\hbar -\frac{i\hbar}{2})} \Psi^\ast_j(z) E_j(u)$$\cr
$$\displaystyle
~~~~~~~~~~~~~~~+ \frac{\mbox{sh}i\pi \eta \hbar}{\pi \eta}
\delta(u-z-\frac{r-l}{2}i\hbar -\frac{i\hbar}{2})
\Psi^\ast_{j-1}(z),$$\cr
$$\displaystyle
E_j(u) \Psi^\ast_{j-1}(z)=\frac{\mbox{sh}\pi \eta(u-z-
\frac{r-j+2}{2}i\hbar - \frac{i\hbar}{2})}{
\mbox{sh}\pi \eta(u-z-
\frac{r-j}{2}i\hbar -\frac{i\hbar}{2})} \Psi^\ast_{j-1}(z)E_j(u),$$\cr
$$E_j(u) \Psi^\ast_{l}(z))=\Psi^\ast_{l}(z) E_j(u),$$~~{\rm otherwise}
\end{array}
\right.\\
& &\displaystyle
 [ F_l(u), \Psi^\ast_j(z)] = \frac{\mbox{sh} i\pi
\eta \hbar}{\pi \eta} \delta_{j,l-1}\delta(u-z-\frac{r-l}{2}i\hbar)
\Psi^\ast_l(z)H^+_l(u+\frac{i\hbar}{4});
\end{eqnarray*}

\item Relations for $\Psi(z)$:

\begin{eqnarray*}
& & \left\{
\begin{array}{l}
$$\displaystyle
\Psi_j(z)H^+_j(u)=\frac{\mbox{sh}\pi \eta(u-z-
\frac{r-j-2}{2}i\hbar -\frac{1}{4}i\hbar)}{
\mbox{sh}\pi \eta(u-z-
\frac{r-j}{2}i\hbar -\frac{1}{4}i\hbar)}H^+_j(u) \Psi_j(z),$$\cr
$$\displaystyle
\Psi_{j-1}(z)H^+_j(u)=\frac{\mbox{sh}\pi \eta(u-z-
\frac{r-j+2}{2}i\hbar -\frac{1}{4}i\hbar)}{
\mbox{sh}\pi \eta(u-z-
\frac{r-j}{2}i\hbar -\frac{1}{4}i\hbar)}H^+_j(u) \Psi_{j-1}(z),$$\cr
$$\displaystyle
\Psi_{l}(z)H^+_j(u)=H^+_j(u) \Psi_{l}(z),$$~~{\rm otherwise}
\end{array}
\right.\\
& & \displaystyle
\left\{
\begin{array}{l}
$$\displaystyle
\Psi_j(z)H^-_j(u)=\frac{\mbox{sh}\pi \eta(u-z-
\frac{r-j-2}{2}i\hbar -\frac{3}{4}i\hbar)}{
\mbox{sh}\pi \eta(u-z-
\frac{r-j}{2}i\hbar -\frac{3}{4}i\hbar)}H^-_j(u) \Psi_j(z),$$\cr
$$\displaystyle
\Psi_{j-1}(z)H^-_j(u)=\frac{\mbox{sh}\pi \eta(u-z-
\frac{r-j+2}{2}i\hbar -\frac{3}{4}i\hbar)}{
\mbox{sh}\pi \eta(u-z-
\frac{r-j}{2}i\hbar -\frac{3}{4}i\hbar)}H^-_j(u) \Psi_{j-1}(z),$$\cr
$$\Psi_{l}(z)H^-_j(u)=H^-_j(u) \Psi_{l}(z),$$~~{\rm otherwise}
\end{array}
\right.\\
& &\left\{
\begin{array}{l}
$$\displaystyle
\Psi_j(z)E_j(u)=\frac{\mbox{sh}\pi \eta(u-z-
\frac{r-j-2}{2}i\hbar -\frac{i\hbar}{2})}{
\mbox{sh}\pi \eta(u-z-
\frac{r-j}{2}i\hbar -\frac{i\hbar}{2})}E_j(u) \Psi_j(z),$$\cr
$$\displaystyle
\Psi_{j-1}(z)E_j(u)=\frac{\mbox{sh}\pi \eta(u-z-
\frac{r-j+2}{2}i\hbar - \frac{i\hbar}{2})}{
\mbox{sh}\pi \eta(u-z-
\frac{r-j}{2}i\hbar -\frac{i\hbar}{2})}E_j(u) \Psi_{j-1}(z)$$\cr
$$\displaystyle
~~~~~~~~~~~~~~~+ \frac{\mbox{sh}i\pi \eta \hbar}{\pi \eta}
\delta(u-z-\frac{r-l}{2}i\hbar - \frac{i\hbar}{2}) \Psi_{j}(z),$$\cr
$$\Psi_{l}(z)E_j(u)=E_j(u) \Psi_{l}(z),$$~~{\rm otherwise}
\end{array}
\right.\\
& & \displaystyle
 [ \Psi_j(z), F_l(u) ] = \frac{\mbox{sh} i \pi \eta \hbar}{
\pi \eta} \delta_{lj} \delta(u-z-\frac{r-l}{2}i\hbar)
H^+_l(u+\frac{i\hbar}{4}) \Psi_{l-1}(z).
\end{eqnarray*}
\end{itemize}

\noindent Similar relations for $q$-affine algebras can be find
in ref.\cite{Dingi}.

\begin{rem}
The intertwining relations are highly sensitive to the form
of the comultiplication used in the definition of intertwining operators.
The relations given above can be obtained only if we use
the comultiplications defined in (\ref{CoB}-\ref{CoE}). For
other form of the comultiplication such as the one used in
\cite{KLP} for $g=sl_2$, not all of the intertwining relations can be
written explicitly.
\end{rem}

Using the bosonic Heisenberg algebra ${\cal H}(\eta)$, one can
in principle obtain bosonic realizations of these intertwining
operators. Then the calculation for the commutation relations
between these intertwining operators and the correlation functions
of such operators will become possible. We leave such tasks
to future studies.

\section{Discussions}

In closing this paper we give the conclusions and some discussions.

We defined the algebra \Alg and its infinite Hopf family for all the simply-laced
Lie algebras $g$. Using the deformed Heisenberg algebra ${\cal H}(\eta)$,
we obtained the level 1 bosonic representation, and then by repeated
use of the comultiplication we get the representations for all
positive integer levels. For $g=sl_{r+1}$, we also gave the simplest
$(r+1)$-dimensional evaluation representation and the intertwining
relations for the level 1 representation and the
$(r+1)$-dimensional evaluation representation.

Clearly, many relevant problems are still left open and among which
we mention several which we would like to solve in future works.

The first problem is: why not non-simply-laced Lie algebras $g$?
Indeed, no reason can be stated {\it a priori} that no analogous
algebras exist for non-simply-laced Lie algebras $g$. However,
for self-consistence we intentionally excluded non-simply-laced
$g$ in our consideration. The reason is that, for such a $g$, the
Cartan matrix is {\it not} symmetric, so that the
Heisenberg algebra ${\cal H}(\eta)$ is not well-defined
(the condition (\ref{sym}) is violated). Probably the way around is to
use the {\it symmetrized} Cartan matrix instead of the Cartan matrix.
Then we can give well-defined Heisenberg algebra ${\cal H}(\eta)$,
but the Serre-like relations (\ref{Serr1}-\ref{Serr2}) are
still not enough to define the algebra \Alg, because there are cases for
$A_{ij}=-2,-3$, etc.

The second problem is the other possible realizations of the algebra
\Alg. For $q$-affine algebras, Yangian doubles and \Al, three
different realizations are know to exist, i.e. the current
realization, Drinfeld generator realization and the
Reshetkhin-Semenov-Tian-Shansky ($RLL$) realization. For our algebra,
it seems important to find the third realization because this
realization has direct connection with the Yang-Baxter relation and
hence is more convenient while considering the possible application
of the algebra in integrable quantum field theories.

As mentioned is the introduction, we postulate that our algebra might
have important application in describing the quantum symmetries
of affine Toda theory, however such applications can be
made possible only if we have identified the $R$-matrix of
our algebra with the quantum $S$-matrix of affine Toda theory.
In this respect, the other form of the comultiplication which
is compatible with $RLL$ relations is also important because under
such a comultiplication the commutation relations
between the intertwining operators would become a set of
Faddeev-Zamolodchikov like algebra which should be explained
as the operator form of the quantum scattering of the corresponding
integrable quantum field theory--the affine Toda theory as
we postulate.

Various considerations on the different choices of domains for
the deformation parameters $\hbar$ and $\eta$ are also
important. On this point the authors of ref.\cite{KLP} have already
listed many problems to whom we whole heartily agree. Besides
the problems listed in there, we are also interested in the
case of $\hbar \rightarrow \infty$, which should correspond
to the case of crystal base for $q$-affine algebras.

Last, we would like to mention the possible connections between
our algebra and the quantum $(\hbar, \xi)$-deformed Virasoro
and $W$-algebras. The $q$- and $\hbar$-deformed Virasoro (and $W$)
Poisson algebras were known to be closely connected to
$q$-affine algebra and Yangian double at the critical level
\cite{FR,DHZ}. The quantum versions of these deformed algebras were also
known to exist and nobody knows to which deformed affine
algebras they correspond. The algebras given in this paper
may be the right candidate to correspond to the quantum
$(\hbar, \xi)$-deformed Virasoro and $W$-algebras.
We point out that algebras correspond to the $(q,p)$-deformed
quantum $W$-algebras also exist, which are
generalizations of the elliptic algebra \Apq to other $g$
with higher rank. We shall present the bosonic representation for
the current realization for such algebras (which we call \Apqg,
representative of yet another example of infinite Hopf family
structure) in the next paper \cite{HZ}.

\end{document}